\documentclass[12pt,oneside]{article}
\usepackage{epsfig}
\usepackage{latexsym}
\usepackage{amsmath}
\usepackage{amssymb}
\usepackage[english]{babel}
\usepackage{xcolor}

\newtheorem{thm}{Theorem}[subsection]
\newtheorem{cor}[thm]{Corollary}
\newtheorem{lem}[thm]{Lemma}
\newtheorem{prop}[thm]{Proposition}

\numberwithin{equation}{section}

\newfont{\sfl}{cmssi12}

\begin{document}

\title{ Quasi-Leontief Utility Functions   on Partially Ordered Sets\, II: Nash Equilibria }

\author{\thanks {University of Perpignan, Department of Economics, 52 avenue Paul Alduy,
66800 Perpignan, France.}  Walter Briec,  QiBin Liang and Charles
Horvath \thanks {University of Perpignan, 52 avenue Paul Alduy,
66800 Perpignan, France.}}  \maketitle \
\\

\bigskip
\textbf{Abstract:} We prove  that, under appropriate conditions, an abstract game with quasi-Leontief payoff functions 
$u_i : \prod_{j=1}^nX_j\to\mathbb{R}$  has a Nash equilibria. When all the payoff functions are globally quasi-Leontief, the existence and the characterization of efficient Nash equilibria  mainly follows from the analysis carried out in part I. When the payoff functions are individually quasi-Leontief functions the matter is somewhat more complicated. We assume that all the strategy spaces are compact topological semilattices, and under appropriate continuity conditions on the payoff functions, we show that there exists an efficient Nash equilibria using the Eilenberg-Montgomery Fixed Point Theorem for acyclic valued upper semicontinuous maps defined on an absolute retract and some non trivial properties of topological semilattices. The map in question is defined on the set of Nash equilibria and its fixed points are exactly the efficient Nash equilibria. 

\newpage
\section{Introduction}\label{secnash}
     Given a family of sets and functions  $S_i\subset X_i$ and $u_i: \prod_{j\in [n]}X_j\to\Lambda$,  a Nash points of the abstract game ${\mathcal G} = \big(u_i, S_i, X_i\big)_{i\in [n]}$ is an element $x^\star$ of the product space such $ \prod_{j\in [n]}S_j$ that, for all $j\in [n]$, 
     $x^\star_j\in\arg\!\max(u_j[x^\star_{-j}]; S_j)$; $\boldsymbol{Nash}\big({\mathcal G}\big)$ denotes the possibly empty set of Nash Points of the abstract game ${\mathcal G} = \big(u_i, S_i, X_i\big)_{i\in [n]}$. In this section we investigate the existence of Nash points in the context of quasi-Leontief functions. If all the payoff functions $u_i$ are individually (respectively globally) quasi-Leontief functions we will say that  ${\mathcal G}$ is an individually (respectively globally) quasi-Leontief game. A quasi-Leontief game is a game which is indifferently either individually or globally quasi-Leontief. Of course, anything that is true of all individually quasi-Leontief games is also true of globally quasi-Leontief games.  
     We will write  ${\mathcal G} = \big(u_i,  X_i\big)_{i\in [n]}$ for  an abstract game for which, for all $i\in [n]$, $S_i = X_i$.

      If each 
     $S_i$ has a largest element $\bar{z}_i$ then the problem of the existence of Nash points is trivial and  
     $\bar{z}$ is a Nash Point. 
     
      \bigskip\noindent We will say that $x^\star = (x_1^\star, \cdots, x_n^\star)$ is an  {\bf efficient  Nash point for player $\boldsymbol{i}$}
       if it is a Nash point and $x_i^\star\in \mathcal{E}(u_i[x_{-i}^\star]_{\mid S_{i}}; S_i)$ and     
       that $x^\star$ is an {\bf efficient Nash point}  if it is efficient for all the players.

     \bigskip
     In section \ref{secnashglob} the strategy spaces are partially ordered spaces, or inf-semilattices, and the payoff functions are globally quasi-Leontief on the product of the strategy spaces; all the results follow   from the previous analysis of quasi-Leontief functions on partially  ordered spaces. Assuming that the constraint sets $S_i$ are comprehensive and bounded above subsets of infsemilattices $X_i$, on can characterize efficient Nash points. 
          
     \medskip In section \ref{secnashindivid} the payoff functions are individually quasi-Leontief and the structure of the strategy spaces is much more restricted, but more classical; the strategy spaces are compact toplogical spaces endowed with an infsemilattice structure for which the inf operation is continuous, as is the case, for example, for compact sub inf-semilattices of $\mathbb{R}^n$. Also, the payoff functions are real valued and continuous. The existence of Nash points in this context does not follow from any of the  previously established results and requires some  topological machinery. We prove the existence of efficient Nash points using the Eilenberg-Montgomery Theorem and some non trivial facts on the topology of inf-semilattices.

     \section{Globally quasi-Leontief games}\label{secnashglob}
    
    Finding a Nash point for a globally quasi-Leontief game with $n$ players  can be reduced to $n$ independent maximization problems. This is content of Proposition \ref{Existnashgql} below     which     settles the problem of the existence of Nash points for globally quasi-Leontief games under the hypothesis that for all $j\in [n]$ 
     $\arg\!\max(u_i; \prod_{j\in [n]}S_j)\neq\emptyset$. Let $S =  \prod_{j\in [n]}S_j$ and for all $i\in [n]$ let $\arg\!\max(u_i; S)_i$ be the projection of 
     $\arg\!\max(u_i; S)$ on $S_i$. 
     
     \begin{prop}\label{Existnashgql} For all globally quasi-Leontief games 
     ${\mathcal G} = \big(u_i, S_i, X_i\big)_{i\in [n]}$ we have 
     \begin{equation*}
     \prod_{i\in [n]}\arg\!\max(u_i; S)_i\, \subset \, \boldsymbol{Nash}\big({\mathcal G}\big).
     \end{equation*}
     \end{prop}
     
     {\it Proof:} If one the sets $\arg\!\max(u_i; S)$ is empty there is nothing to prove. For all $i\in [n]$ pick an element $z^{[i]}$ of $\arg\!\max(u_i; S)$ and let $x^\star_i = z^{[i]}_i$. We will show that $x^\star$ is a Nash point. 
 For for all $x_i\in S_i$ we have  $u_i(z^{[i]})\geqslant u_i(x^\star_{-i}; x_i)$, that is  
     $u_i[z^{[i]}_{-i}](z_i^{[i]})\geqslant u_i[x^{\star}_{-i}](x_i)$; and trivially we have $ u_i[x^{\star}_{-i}](x_i)\geqslant  u_i[x^{\star}_{-i}](x_i)$. From the remarks following Proposition $3.1.2$ of Part I, 
      we have 
     \begin{equation}\label{eqnasgql1}u_i[z^{[i]}_{-i}]^\sharp\big( u_i[x^{\star}_{-i}](x_i)\big) = 
     u_i[x^{\star}_{-i}]^\sharp\big( u_i[x^{\star}_{-i}](x_i)\big) = 
     u_{i, i}^\sharp\big( u_i[x^{\star}_{-i}](x_i)\big)
     \end{equation}
     where $u_{i, i}^\sharp$ is the $i$-th coordinate of $u_{i}^\sharp$. Again, from $u_i[z^{[i]}_{-i}](z_i^{[i]})\geqslant u_i[x^{\star}_{-i}](x_i)$ we have 
      \begin{equation}\label{eqnasgql2}
     z_i^{[i]}\geqslant  u_i[z^{[i]}_{-i}]^\sharp\big(u_i[x^{\star}_{-i}](x_i)\big).
     \end{equation}
     From (\ref{eqnasgql1}) and (\ref{eqnasgql2}) we get 
     $z_i^{[i]}\geqslant u_i[x^{\star}_{-i}]^\sharp\big( u_i[x^{\star}_{-i}](x_i)\big)$ or $u_i[x^{\star}_{-i}](z_i^{[i]})\geqslant u_i[x^{\star}_{-i}](x_i)$ that is 
     $u_i(x^\star)\geqslant u_i(x^{\star}_{-i}; x_i)$. \hfill$\Box$

     \bigskip\noindent
     The next results shows that, under appropriate but by now familiar conditions, a globally quasi-Leontief game has a Nash point which is also a maximal element of the strategy space. 
     
     \begin{prop}\label{Enashmax} 
      Let ${\mathcal G} = \big(u_i, S_i, X_i\big)_{i\in [n]}$ be a globally quasi-Leontief game such that, for all $i\in [n]$, $S_i$ is a non empty comprehensive (CUC) subset of $X_i$ with an  upper bound  $\bar{x}_i\in X_i$. Then 
       \begin{equation*}
    \boldsymbol{Max}\big( \prod_{i\in [n]}S_i\big)\, \cap \, \boldsymbol{Nash}\big({\mathcal G}\big) \, \neq\, \emptyset.
     \end{equation*}
     \end{prop}
     
      {\it Proof:} The set $S = \prod_{i\in [n]}$ is clearly non empty, comprehensive and bounded above; one easily shows that it is also (CUC) since a projection of a chain in the product space is a chain. By Theorem $2.4.3$ of Part I,  
      $\boldsymbol{Max}(S)\cap \arg\!\max(u_i; S)\neq\emptyset$, for all 
      $i\in [n]$; in the proof of Proposition \ref{Existnashgql} we take 
      $z^{[i]}$ in  $\boldsymbol{Max}(S)\cap \arg\!\max(u_i; S)$. Let us see that $x^\star\in \boldsymbol{Max}(S)$. Let $x\in S$ such that $x\geqslant x^\star$; from $x_i\geqslant z_i^{[i]}$ we have $(z_{-i}^{[i]}; x_i)\geqslant z^{[i]}$ and from $(z_{-i}^{[i]}; x_i)\in S$ we obtain $(z_{-i}^{[i]}; x_i) = z^{[i]}$ and consequently $x_i = z_i^{[i]}$, and since $i$ was arbitrary, $x = x^\star$. \hfill$\Box$

      \bigskip\noindent
  \subsection{Efficient Nash points of globally quasi-Leontief games.} 
     In this section, each $X_i$ is an inf-semilattice and  each  $S_i$ is a nonempty comprehensive subset of $X_i$ with an upper bound $\bar{x}_i\in X_i$. By Proposition $2.3.5$ of Part I,    there are quasi-Leontief functions $u_{i, j} : S_i\to\Lambda$, $i, j\in [n]$ such that, for all $x\in\prod_{i\in [n]}S_i$, $u_i(x) = \min_{j\in [n]}u_{i, j}(x_j)$. For all $i\in [n]$ let 
     $\tilde{u}_i(x_{-i}) = \min_{j\neq i}u_{i, j}(x_j)$; it is a globally quasi-Leontief function defined on 
     $\prod_{j\neq i}X_j$. We have $u_i(x_1, \cdots, x_n) = \min\{u_{i, i}(x_i), \tilde{u}_i(x_{-i})\}$.
     
     \medskip\noindent
     A point $x^\star = (x^\star_1, \cdots, x^\star_n)$ is a Nash point if and only if, for all $i\in [n]$, the following condition holds: 
     
     \centerline{$(N_i)$\quad$\forall x_i\in X_i$\, $u_{i, i}(x^\star_i)\, \geqslant\,  \min\{u_{i, i}(x_i), \tilde{u}_i(x^\star_{-i})\} $.}
     
     \medskip\noindent
     If $x_i^\star\in \arg\!\max(u_{i, i}; S_i)$ then $(N_i)$ holds. If $x^\star$ is a Nash point such that $x_i^\star\not\in \arg\!\max(u_{i, i}; S_i)$ then there exists $x_i\in X_i$ such that $u_{i, i}(x_i) > u_{i, i}(x^\star_i)$ and threfore, by $(N_i)$, $u_{i, i}(x^\star_i)\, \geqslant\,   \tilde{u}_i(x^\star_{-i})$; and this inequality obviously imply $(N_i)$. \\ If the inequality  $u_{i, i}(x^\star_i)\, \geqslant\,   \tilde{u}_i(x^\star_{-i})$ does not hold then $(N_i)$ implies that $x_i^\star\in \arg\!\max(u_{i, i}; S_i)$.\\
     If $x_i^\star\in \arg\!\max(u_{i, i}; S_i)$ then $u_{i, i}^\circ(x_i^\star)\in \arg\!\max(u_{i, i}; S_i)$, since $S_i$ is comprehensive. \\
     If $u_{i, i}(x^\star_i)\, \geqslant\,  \min_{j\neq i}\{u_{i, j}(x^\star_j)\}$ then, from $u_{i, i}(x^\star_i) = u_{i, i}\big(u_{i, i}^\circ(x^\star_i)\big)$,\\ $x^\star_j\geqslant u_{j, j}^\circ(x^\star_j)$ and the fact that $u_{i, j}$ is isotone we have\\ 
     $ u_{i, i}\big(u_{i, i}^\circ(x^\star_i)\big)\, \geqslant\,  \min_{j\neq i}\{u_{i, j}(x^\star_j)\} \geqslant\,  \min_{j\neq i}\big\{u_{i, j}\big(u_{j, j}^\circ(x^\star_j)\big)\big\}$.
     
     \begin{prop}\label{carnashglobleon} If, for all $i\in [n]$, $S_i$ is a comprehensive subset of the inf-semilattice $X_i$ and if,  for all $i\in [n]$, 
     $u_i(x_1, \cdots, x_n) = \min_{i\in[n]}u_{i, j}(x_j)$, where $u_{i, j}: X_j\to\Lambda$ is a quasi-Leontief function, then a point $x^\star\in\prod_{i\in [n]}X_i$ is a Nash point if and only if, for all $i\in [n]$, 
     \begin{equation}\label{N1}
     x_i^\star\in \arg\!\max(u_{i, i}; S_i)
     \end{equation} or
     \begin{equation}\label{N2}u_{i, i}(x^\star_i)\, \geqslant\,   \tilde{u}_i(x^\star_{-i})
\end{equation}
Furthermore, if $(x_1^\star, \cdots, x_n^\star)$ is a Nash point then so is 
$(u_{1, 1}^\circ(x_1^\star), \cdots, u_{n, n}^\circ (x_n^\star) )$ and, for all 
$i\in [n]$, $u_i((x_1^\star, \cdots, x_n^\star) \geqslant u_i(u_{1, 1}^\circ(x_1^\star), \cdots, u_{n, n}^\circ (x_n^\star) )$.

\end{prop}

\noindent The meaning of Proposition \ref{carnashglobleon} is that the search of a Nash point for a globally quasi-Leontief game on  comprehensive and bounded strategy spaces $S_i$ can always be reduced to $n$ independent and decoupled maximization problems -- 
$x_i^\star\in \arg\!\max(u_{i, i}; S_i)$ -- and that one can choose 
$x_i^\star\in\mathcal{E}(u_{i, i}; S_i)$.

\bigskip Recall that $x^\star$ is an efficient Nash point for player $i$ if 
$x_i^\star\in\mathcal{E}(u_i[x^\star_{-i}]_{\mid S_{i}}, S_i)$ that is, if $\big [x_i\in S_i$ and $u_i[x^\star_{-i}](x_i)\geqslant u_i(x^\star)\big]$ implies 
$x_i \geqslant x_i^\star$,  which is equivalent to 
\begin{equation}\label{condeffnash}
\big[ x_i\in S_i \hbox{ and } u_{i, i}(x_i) \geqslant \min\{u_{i, i}(x^\star_i), \tilde{u}_i(x^\star_{-i}) \}\big] \hbox{ implies } x_i\geqslant x_i^\star
\end{equation}

\noindent   $(a1)$ If $\tilde{u}_i(x^\star_{-i}) \geqslant u_{i, i}(x^\star_i)$ then (\ref{condeffnash}) becomes: $\big[ x_i\in S_i \hbox{ and } u_{i, i}(x_i) \geqslant u_{i, i}(x^\star_i)]\\ \hbox{ implies } x_i\geqslant x_i^\star$
.\\ 
If $x_i^\star\not\in\mathcal{E}(u_{i, i}; S_i)$ then $u_{i, i}^\circ(x_i^\star)\in S_i$ and 
$u_{i, i}^\circ(x_i^\star)\not\in\uparrow\!(x^\star_i)$, since $u_{i, i}^\circ(x_i^\star)\neq x_i^\star$, and 
$u_{i, i}\big(u_{i, i}^\circ(x_i^\star)\big) = u_{i, i}(x^\star_i) =  \min\{u_{i, i}(x^\star_i), \tilde{u}_i(x^\star_{-i}) \}$ and therefore (\ref{condeffnash}) does not hold.\\
If $x_i^\star\in\mathcal{E}(u_{i, i}; S_i)$ then (\ref{condeffnash}) holds. 

\medskip\noindent
{\it A Nash point $x^\star$ such that $\tilde{u}_i(x^\star_{-i}) \geqslant u_{i, i}(x^\star_i)$ is efficient for player $i$ if and only if  $x_i^\star\in\mathcal{E}(u_{i, i}; S_i)$.  We also have $u_i(x^\star) = u_{i, i}(x^\star_i)$. }

\medskip\noindent
$(a2)$ If $u_{i, i}(x^\star_i) > \tilde{u}_i(x^\star_{-i})$ then (\ref{condeffnash}) becomes: 
$\big[ x_i\in S_i \hbox{ and } u_{i, i}(x_i) \geqslant \tilde{u}_i(x^\star_{-i}) \big]$ implies  $x_i\geqslant x_i^\star$.\\ If $x_i^\star\not\in\mathcal{E}(u_{i, i}; S_i)$ then $u_{i, i}\big(u_{i, i}^\circ(x_i^\star)\big) = u_{i, i}(x^\star_i) > \tilde{u}_i(x^\star_{-i})$, $u_{i, i}^\circ(x_i^\star)\in S_i$ and 
$u_{i, i}^\circ(x_i^\star)\not\in\uparrow\!(x^\star_i)$ and therefore (\ref{condeffnash}) does not hold.\\
If $x_i^\star\in\mathcal{E}(u_{i, i}; S_i)$ and (\ref{condeffnash}) does not hold then there exists 
       $x_i\in S_i$ such that $x_i\not\in\uparrow(x_i^\star)$ and 
       $u_{i, i}(x_i)  \geqslant  \tilde{u}_i(x^\star_{-i})
       $.\\
       We cannot have $u_{i, i}(x_i)\geqslant u_{i, i}(x_i^\star) $ since this would imply $x_i\geqslant x_i^\star$; therefore 
        $u_{i, i}(x_i^\star) >  u_{i, i}(x_i)$. We have shown that          \begin{equation}\label{cnsnoteffnash}
        \exists x_i\in S_i \hbox{ such that } 
        u_{i, i}(x_i^\star) >  u_{i, i}(x_i)  \geqslant \tilde{u}_i(x^\star_{-i})
       .
        \end{equation}
        Reciprocally, if (\ref{cnsnoteffnash}) holds we cannot have $x_i\geqslant x_i^\star$, since $u_{i, i}$ isotone,  and therefore (\ref{condeffnash}) does not hold.
        
        \medskip\noindent{\it A Nash point $x^\star$ such that $u_{i, i}(x^\star_i) > \tilde{u}_i(x^\star_{-i})$ is efficient for player $i$ if and only if $x_i^\star\in\mathcal{E}(u_{i, i}; S_i)$ and, for all 
        $x_i\in S_i$, either $u_{i, i}(x_i)\geqslant u_{i, i}(x_i^\star)$ or 
        $ \tilde{u}_i(x^\star_{-i})\geqslant u_{i, i}(x_i)$. We also have 
        $u_i(x^\star) =  \tilde{u}_i(x^\star_{-i})$.    }
        
        \medskip\noindent Now, let us assume that $X_i$ is a topological semilattice and that $S_i$ is a connected subset of $X_i$. If $u_{i, i}$ is continuous on $S_i$ and if $u_{i, i}(x^\star_i) > \tilde{u}_i(x^\star_{-i})$ 	and $ \tilde{u}_i(x^\star_{-i})\geqslant u_{i, i}(x_i)$ for at least one $x_i$ in $S_i$ then there must exists $x_i^\prime\in S_i$ such that 
        $u_{i, i}(x^\star_i) > u_{i, i}(x^\prime_i) > \tilde{u}_i(x^\star_{-i})$. 
        
         \medskip\noindent{\it If $X_i$ is a topological semilattice and that $S_i$ is a connected subset of $X_i$, a Nash point $x^\star$ such that $u_{i, i}(x^\star_i) > \tilde{u}_i(x^\star_{-i})$ is efficient for player $i$ if and only if $x_i^\star\in\mathcal{E}(u_{i, i}; S_i)$ and, for all 
        $x_i\in S_i$,  $u_{i, i}(x_i)\geqslant u_{i, i}(x_i^\star)$, that is \\$x_i^\star\in  \arg\!\min(u_{i, i}; S_i)$.   Therefore, a Nash point $x^\star$ such that $x_i^\star\in  \arg\!\max(u_{i, i}; S_i)$ can not be efficient for player $i$, unless $u_{i, i}$ is constant.  }
        
        \medskip\noindent What is the meaning of the inequality $u_{i, i}(x^\star_i) > \tilde{u}_i(x^\star_{-i})$?   Assume that all the strategy  spaces $X_j$ are topological spaces and that all the functions $u_{i, j}$, $j\in [n]$, are continuous.  One can find, for all $j\in [n]$, a neighbourhood $V^{[i]}_j$ of $x^\star_j$ in $S_j$ such that, for all $x\in V^{[i]} =   \prod_{j\in [n]}V^{[i]}_j$, 
  $u_{i, i}(x_i) >  \tilde{u}_i(x_{-i})
       $ and therefore,  for all 
       $x\in V^{[i]}$,  
       \begin{equation}\label{tildeu}
       u_i(x) =  \tilde{u}_{-i}(x_{-i})
       \end{equation}
       that is:\\
       {\it  if $u_{i, i}(x^\star_i) > \tilde{u}_i(x^\star_{-i})$ then, in a  neighbourhood  $V^{[i]}$ of  $x^\star$, the payoff function $u_i$ of player $i$ is a function of the strategies of the remaining players and not of its own strategy.}

       \medskip\noindent
        If $x^\star$ is such that, for all $i\in [n]$, $u_{i, i}(x_i^\star) > \tilde{u}_{-i}(x^\star_{-i}) $ then there exists a neighbourhood $W$ of $x^\star$ in $\prod_{i\in [n]}S_i$ (for example the intersection of all the  $V^{[j]}$ above), on which, for all $i\in [n]$, the payoff function of player $i$ depends only on the strategies of the remaining players; more precisely,  for all $i\in [n]$ and all $x\in W$, 
       \begin{equation}\label{totineff}
       u_i(x) =  \tilde{u}_{-i}(x_{-i}).
       \end{equation}

     \bigskip\noindent
     Let us look at an example with two players. 
     
      \medskip\noindent
     $S_1 = [0, 2]$, $S_2 = [0, 2]$,  $
u_{1,1}(x_1) = \left\{
\begin{array}{lcl}
 2x_1 & \hbox{if}  &   0\leqslant x_1\leq 1\\
 2 &  \hbox{if} & 1\leqslant x_1 \end{array}
\right.
$,   $u_{1, 2}(x_2) = \displaystyle{\frac{x_2}{2}}$

\medskip\noindent and, $u_{2, 1}(x_1) = \displaystyle{\frac{x_1}{2}}$, 
$u_{2,2}(x_2) = \left\{
\begin{array}{lcl}
 2x_2 & \hbox{if}  &   0\leqslant x_2\leq 1\\
 2 &  \hbox{if} & 1\leqslant x_2 \end{array}
\right.
$.\\
Take  $x_1^\star = 1$, $x_2^\star = 1$; we then have  $u_{1,1}(x_1^\star) = u_{2,2}(x_2^\star)$, $u_{1,2}(x_2^\star) = u_{2,1}(x_1^\star) = \displaystyle{\frac{1}{2}}$ and therefore 
$\left\{
\begin{array}{lll}
 u_1(x_1^\star, x_2^\star) & =    &  u_{1,2}(x_2^\star) =  \displaystyle{\frac{1}{2}}\\
 & & \\
 u_2(x_1^\star, x_2^\star) &  = &  u_{2,1}(x_1^\star) = \displaystyle{\frac{1}{2}} \end{array}
\right.$ 
\medskip\noindent and also $u_1(3/4, x_2^\star) = u_1(x_1^\star, x_2^\star)$ with 
$x_1^\star > 3/4$ and similarly for $u_2$. In conclusion at the Nash point 
$(x_1^\star, x_2^\star) = (1, 1)$ the payoff of player $1$ depends only on the strategy of player $2$ and not on its own strategy and similarly for player $2$. As one can see, $(1, 1)$ is not an efficient Nash point. But, also, 
$(2, 2)$ is a Nash point since it is the largest element of $S_1\times S_2$. 
There is an efficient Nash point: $(0, 0)$. It gives the smallest possible payoff to both players. So, a Nash point that gives the largest possible payoff to both players is not efficient and the Nash point that gives the smallest possible payoff to both players is efficient.   

\bigskip\noindent
The problem of the existence of efficient Nash points for individually - and therefore globally - quasi-Leontief games will be treated in the next section.

\section{Indivually quasi-Leontief games}\label{secnashindivid}
Let us start by  defining  some of the  concepts that will be needed to state the theorem from which the existence of Nash points for  indivually quasi-Leontief games will be deduced. 

\medskip\noindent
We will assume  that $\Lambda = \mathbb{R}$ and that all the strategy spaces $X_i$ are  {\bf topological inf-semilattices} that is: 

\medskip\noindent
{\it $X_i$ is an inf-semilattice endowed with a topology for which the inf operation $\wedge : X_i\times X_i\to X_i$ is continuous.} 

\medskip\noindent We recall that a subspace $Z$ of a topological space $X$ is path connected if, for all pair $(z_0, z_1)\in Z\times Z$, there exists a continuous map $\alpha : [0, 1]\to  Z$ such that $\alpha(0) = z_0$ and 
$\alpha(1) = z_1$.

\medskip\noindent A subset $C$ of an inf-semilattice $X$ is {\bf inf-convex} if, for all  $x_0, x_1\in X$, 
$[x_0\wedge x_1, x_0]\subset X$.

\medskip\noindent An easy induction shows that if $C$ is an inf-convex subset of an inf-semilattice $X$ then, for all finite an non empty subset 
$S\subset C$, $\cup_{x\in S}[\wedge S, x]\subset C$. 

\medskip\noindent Also, a subset $C$ of an inf-semilattice $X$ is inf-convex if and only if the following two conditions hold: 

\medskip\noindent$(1)$ $S$ is a sub-semilattice of $X$; that is, for all $x$ and $y$ in $S$, $x\wedge y\in S$\quad and 

\medskip\noindent$(2)$ for all element $(x, y)\in C\times C$ such that 
$x\leqslant y$ one has $[x, y]\subset C$.

\begin{thm}\label{luotheo} Let ${\mathcal G} = \big(u_i, X_i\big)_{i\in [n]}$ be an abstract game such that: 

\medskip\noindent
$(1)$ the strategy spaces $X_i$ are compact metrizable inf-semilattices with path-connected intervals; 

\medskip\noindent
$(2)$ the payoff functions $u_i : \prod_{j\in [n]}X_j\to\mathbb{R}$ are continuous and such that, for all $x_{-i}\in\prod_{j\neq i}X_j$ and all $t\in\mathbb{R}$, the set 
$\{z_i\in X_i: u_i[x_{-i}](z_i) > t\}$ is inf-convex;

\medskip\noindent Then $\boldsymbol{Nash}\big({\mathcal G}\big) \neq\emptyset$.
\end{thm}

Theorem \ref{luotheo} is a simplified version of  Theorem 4.1 of Luo \cite{luo}.

 \bigskip
 We know from Lemma 2.3.1 of Part I that for a quasi-Leontief $u: X\to\mathbb{R}$ function defined on an inf-semilattice $X$ one always has 
 $u(x_1\wedge x_2) = \min\{u(x_1), u(x_2)\}$ and this implies that, for all 
 $t\in\mathbb{R}$, $\{x\in X: u(x) > t\}$ is inf-convex. From Theorem \ref{luotheo} we have the existence of Nash points for abstract 
 individually quasi-Leontieff games. 

\begin{thm}\label{existsnashindi} If ${\mathcal G} = \big(u_i, X_i\big)_{i\in [n]}$ be an abstract  individually  quasi-Leontieff game such that: 

\medskip\noindent
$(1)$ the strategy spaces $X_i$ are compact metrizable inf-semilattices with path-connected intervals; 

\medskip\noindent
$(2)$ the payoff functions $u_i : \prod_{j\in [n]}X_j\to\mathbb{R}$ are continuous.

\medskip\noindent Then $\boldsymbol{Nash}\big({\mathcal G}\big)\neq\emptyset$.

\end{thm}

\bigskip\noindent
In $\mathbb{R}^m$ with the partial order associated to the positive cone\,  
$\mathbb{R}^m_+$ order intervals are path-connected therefore, an inf-convex subset of  $\mathbb{R}^m$ is an inf-semilattice with path-connected intervals.

\begin{cor}
\label{existsnashindiRm} Let ${\mathcal G} = \big(u_i, X_i\big)_{i\in [n]}$ be an abstract individually  quasi-Leontieff game such that: 

\medskip\noindent
$(1)$ for all $i\in [n]$, the strategy space $X_i$ is a compact inf-convex subsets of 
$\mathbb{R}^{m_i}$; 

\medskip\noindent
$(2)$ the payoff functions $u_i : \prod_{j\in [n]}X_j\to\mathbb{R}$ are continuous.

\medskip\noindent Then $\boldsymbol{Nash}\big({\mathcal G}\big)\neq\emptyset$. 
\end{cor}

\bigskip One can easily see that an arbitrary intersection of inf-convex subsets is inf-convex and that the union $\cup\mathcal{C}$ of a family of inf-convex subsets is inf-convex if, for all $C, C^\prime\in \mathcal{C}$ there exists $C^{\prime\prime}\in\mathcal{C}$ such that $C\cup C^\prime\subset C^{\prime\prime}$; consequently, the following statements are equivalent: 

\medskip\noindent$(a)$ for all $t\in\mathbb{R}$, $\{x\in X: u(x) > t\}$ is 
inf-convex;  

\medskip\noindent$(b)$ for all $t\in\mathbb{R}$, $\{x\in X: u(x) \geqslant t\}$ is 
inf-convex.  

\bigskip\noindent
Let us write $[\![x_1, x_2]\!]$ for $[x_1\wedge x_2, x_1]\cup [x_1\wedge x_2, x_2]$. 
One can now easily see that $C$ is inf-convex if and only if, for all $x_1, x_2\in C$, $[\![x_1, x_2]\!]\subset C$ and that 
condition $(a)$ above is equivalent to
\begin{equation}\label{infquasiconcave}
\forall x_1, x_2\in X\quad\inf_{x\in [\![x_1, x_2]\!]}\!\!u(x)\geqslant\min\{u(x_1), u(x_2)\}.
\end{equation}

\medskip\noindent Luo's Theorem applies to a much larger class than the class of individually quasi-Leontief functions; it only requires the payoff functions to be ``inf-quasiconvex'' in each variable. But Luo's Theorem, which is derived from a Browder-Fan fixed point theorem for topological inf-semilattices, does not say anything about the existence of efficient Nash points. We will show that efficient Nash points always exists but this requires
a Kakutani like fixed point theorem in topological inf-semilattices which can not be, at least as far as we know, established from the single assumption that intervalls are path connected.

\subsection{Efficient Nash points for individually quasi-Leontief games}\label{effnashsubsec}
We assume that the strategy spaces are compact topological inf-semilattices which are metrizable and that the payoff functions 
$u_i : \prod_{j\in [n]}X_j\to\mathbb{R}$ are continuous and individually quasi-Leontieff. We want to show that there exists an efficient Nash point 
$x^\star$, that is, $x^\star\in \boldsymbol{Nash}\big({\mathcal G}\big) $ and  $x^\star\in\prod_{i\in[n]}\mathcal{E}(u_i[x_{-i}^\star], X_i)$. In other words, we want to show that the map $\mathbb{P}_{\mathcal{G}}$ defined by $\mathbb{P}_{\mathcal{G}}(x) = \prod_{i\in[n]}\mathcal{E}(u_i[x_{-i}], X_i)$ has a fixed point in $\boldsymbol{Nash}\big({\mathcal G}\big) $. 

\medskip\noindent
For all $x\in\prod_{i\in [n]}X_i$ let 
\begin{equation*}
\mathbb{E}(x) = \prod_{i\in [n]}\arg\!\max(u_i[x_{-i}], \mathcal{E}(u_i[x_{-i}], X_i)).
\end{equation*} 

\begin{lem}\label{xdansEx}
A point $x^\star\in\prod_{i\in [n]}X_i$ is an efficient Nash point if and only if 
$x^\star\in\mathbb{E}(x^\star)$
\end{lem}
{\it Proof:} If $x^\star\in\mathbb{E}(x^\star)$ then, for all $i\in[n]$ and for all $z_i\in\mathcal{E}(u_i[x^\star_{-i}], X_i)$, $u_i(x^\star)\geqslant u_i[x^\star_{-i}](z_i)$. If $x_i$ is an arbitrary element of $X_i$ then 
$u_i[x^\star_{-i}]^\circ(x_i)\in \mathcal{E}(u_i[x^\star_{-i}], X_i)$ and 
$ u_i[x^\star_{-i}](x_i) = u_i[x^\star_{-i}]\big(u_i[x^\star_{-i}]^\circ(x_i)\big)$. This shows that $x^\star$ is a Nash point. From $\mathbb{E}(x^\star)\subset\mathbb{P}_{\mathcal{G}}(x^\star)$, $x^\star$ is an efficient Nash point. 

\medskip\noindent
An efficient Nash points $x^\star$ belongs to 
$\prod_{i\in [n]}\arg\!\max(u_i[x^\star_{-i}],  X_i)$ and to 
$\prod_{i\in [n]} \mathcal{E}(u_i[x_{-i}], X_i)$ and therefore to 
$\mathbb{E}(x^\star)$.\hfill$\Box$

\bigskip The set $\mathbb{E}(x)$ does not have to be an inf-convex subset  of the product space since $x_1\leqslant z \leqslant x_2$ with, for all $j\in [n]$ and $i\in\{1, 2\}$, 
$x_{i, j}\in \mathcal{E}(u_i[x_{-j}], X_j)$ does not imply $z_{j}\in \mathcal{E}(u_j[x_{-j}], X_j)$ - being between two efficient points does not imply efficiency - all we have is $u_j[x_{-j}]^\circ(z_{j})\in \mathcal{E}(u_j[x_{-j}], X_j)$. And this is the cause of some  complications. The structure of  $\mathbb{E}(x)$ is the subject matter of the following short sequence of lemmas.

\begin{lem}\label{pasemptyE} For all $x\in\prod_{i\in [n]}X_i$ the set 
$\mathbb{E}(x)$ is  not empty. 

\end{lem}

{\it Proof:} Since $u_i[x_{-i}]: X_i\to\mathbb{R}$ is continuous and $X_i$ is compact $\arg\!\max(u_i[x_{-i}], X_i)\neq\emptyset$; if $z_i\in \arg\max(u_i[x_{-i}], X_i)$ then $u_i[x_{-i}]^\circ(z_i)\in\arg\!\max(u_i[x_{-i}], X_i)\cap\mathcal{E}(u_i[x_{-i}], X_i)$. \hfill$\Box$

\begin{lem}\label{expath} Assume that the strategy spaces are all compact  inf-semilattices with path-connected intervals. If, for all $x\in\prod_{i\in[n]}X_i$ and all $i\in [n]$,\\ $u_i[x_{-i}]^\circ : X_i\to X_i$ is continuous then, for all $x\in\prod_{i\in[n]}X_i$, the set $\mathbb{E}(x)$ is a  topological inf-semilattice with path connected intervals. Furthermore, $\mathbb{E}(x)$ is compact and it has a smallest and a largest element.
\end{lem}
{\it Proof:} From Lemma 2.3.2 of Part I, $\mathcal{E}(u_i[x_{-i}], X_i)$ is a sub-semilattice of $X_i$; if $S_i$ is a sub-semilattice of $X_i$ then, from $u_i[x_{-i}](x_1\wedge x_2) = \min\{u_i[x_{-i}](x_1), u_i[x_{-i}](x_2)\}$, 
$\arg\!\max(u_i[x_{-i}], S_i)$ is also a sub-semilattice of $X_i$ and therefore\\ $\arg\!\max(u_i[x_{-i}], \mathcal{E}(u_i[x_{-i}], X_i)$ is a sub-semilattice of $X_i$. 

\medskip\noindent
 Since the topology and the inf-operation on $\mathbb{E}(x)$ are those induced from $\prod_{i\in [n]}X_i$ the inf-operation restricted to $\mathbb{E}(x)$ is continuous. We have shown that the  product space $\mathbb{E}(x)$ is a sub-semilattice of the product $\prod_{i\in [n]}X_i$.

\medskip\noindent
Let $x_0$ and $x_1$ be two elements of $\mathbb{E}(x)$ such that 
$x_0\leqslant x_1$. The order interval in  $\mathbb{E}(x)$ is 
$\{z\in \mathbb{E}(x): x_0\leqslant z\leqslant x_1\}$ that is 
$[x_0, x_1]\, \cap \, \mathbb{E}(x)$ where $[x_0, x_1]$ is the order interval in $\prod_{i\in [n]}X_i$. We have to show that there exists a continuous path 
$\alpha: [0, 1]\to \prod_{i\in [n]}X_i$ such that $\alpha(i) = x_i$ for 
$i\in\{0, 1\}$ and, for all $t\in [0, 1]$, $\alpha(t)\in [x_0, x_1]\, \cap \, \mathbb{E}(x)$.\\ 
Choose a continuous path $\gamma_i : [0, 1]\to [x_{0, i}, x_{1, i}]$ and let  
$\alpha_i(t) = u_i[x_{-i}]^\circ(\gamma_i(t))$; we have $\alpha_i(t)\in  \mathcal{E}(u_i[x_{-i}], X_i)$.\\
Since $x_{0, i}$ and $x_{1, i}$ are both in $\arg\!\max(u_i[x_{-i}], X_i)$ with 
$x_{0, i},\leqslant x_{1, i}$ and since $u_i[x_{-i}]$ is isotone we have 
$[x_{0, i}, x_{1, i}]\subset \arg\!\max(u_i[x_{-i}], X_i)$ and consequently 
$\gamma_i(t)\in\arg\!\max(u_i[x_{-i}], X_i)$. From 
$u_i[x_{-i}](\gamma_i(t)) = u_i[x_{-i}]\big(u_i[x_{-i}]^\circ(\gamma_i(t))\big)$ we  have $\alpha_i(t)\in\arg\!\max(u_i[x_{-i}], X_i)$. \\In conclusion, 
$\alpha_i(t)\in \arg\!\max(u_i[x_{-i}], \mathcal{E}(u_i[x_{-i}], X_i) )$. 

\medskip\noindent
 Since $u_i[x_{-i}]^\circ$ is continuous and $X_i$ is compact the set of fixed points of $u_i[x_{-i}]^\circ$, that is $ \mathcal{E}(u_i[x_{-i}], X_i)$, is compact; from the continuity of $u_i[x_{-i}]$ follows the compactness of 
 $\arg\!\max(u_i[x_{-i}], \mathcal{E}(u_i[x_{-i}], X_i) )$. This shows that 
 $\mathbb{E}(x)$ is compact. 
 
 \medskip\noindent We have already seen that a compact inf-semilattice has a smallest element, let $\boldsymbol{m}(x)$ be the smallest element of 
 $\mathbb{E}(x)$. But $\arg\!\max(u_i[x_{-i}], \mathcal{E}(u_i[x_{-i}], X_i) )$ is totally ordered since it is a subset of the totally ordered set $ \mathcal{E}(u_i[x_{-i}], X_i)$; by compactness it has a largest element 
 $\boldsymbol{M}_i(x)$;  $\boldsymbol{M}(x)$ is the largest element of $\mathbb{E}(x)$.
\hfill$\Box$

\begin{lem}\label{scsE} Assume that the strategy spaces are all compact metizable inf-semilattices. If $\mathbb{P}_{\mathcal{G}}$ is upper semicontinuous then $\mathbb{E}$ is upper semicontinuous. 
\end{lem}
{\it Proof:} By compactness it is sufficient to show that the graph of 
$\mathbb{E}$ is a closed subset of $\prod_{i\in[n]}X_i\, \times\, \prod_{i\in[n]}X_i$. We have $(x, y)\in\mathbb{E}$ if and only if, for all $i\in [n]$, 
$y_i\in \mathcal{E}(u_i[x_{-i}], X_i)\cap\arg\!\max(u_i[x_{-i}], X_i)$.  Let 
$(x_{\boldsymbol{m}}, y_{\boldsymbol{m}})_{\boldsymbol{m}\in\mathbb{N}}$ be a sequence of elements of the graph of $\mathbb{E}$ that converges to a point $(\bar{x}, \bar{y})$.  
From $y_{\boldsymbol{m}, i}\in\arg\!\max(u_i[x_{\boldsymbol{m},-i}], X_i)$ and the continuity of $u_i$ we have 
$\bar{y}_{i}\in\arg\!\max(u_i[\bar{x}_{-i}], X_i)$.

\medskip\noindent
The sequence $(x_{\boldsymbol{m}}, y_{\boldsymbol{m}})_{\boldsymbol{m}\in\mathbb{N}}$ also belongs to the graph of $\mathbb{P}_{\mathcal{G}}$, which  is upper semicontinuous; therefore  $\bar{y}_i\in \mathcal{E}(u_i[\bar{x}_{-i}], X_i)$. \hfill$\Box$

\begin{lem}\label{scsP} Assume that the strategy spaces are all compact metizable inf-semilattices, and that: 

\medskip\noindent
$(1)$ for all $x\in\prod_{i\in[n]}X_i$ and all $i\in [n]$, $u_i[x_{-i}]^\circ $ is continuous; 

\medskip\noindent
$(2)$ for all $i\in [n]$ and for all convergent sequence 
$(x_{\boldsymbol{m}})_{\boldsymbol{m}\in\mathbb{N}}$ of points of the product space  
$\prod_{i\in[n]}X_i$ with limit $\bar{x}$, the sequence of functions 
$(u_i[x_{\boldsymbol{m}, -i}]^\circ)_{\boldsymbol{m}\in\mathbb{N}}$ converges pointwise to $u_i[\bar{x}_{ -i}]^\circ$.

\medskip\noindent Then $\mathbb{P}_{\mathcal{G}}$  is upper semicontinuous.

\end{lem}
{\it Proof:} Assume that the sequence $(x_{\boldsymbol{m}}, y_{\boldsymbol{m}})_{\boldsymbol{m}\in\mathbb{N}}$ converges to 
$(\bar{x}, \bar{y})$ and that $y_{\boldsymbol{m}}\in\mathbb{P}_{\mathcal{G}}(x_{\boldsymbol{m}})$. We have to see that, for all 
$i\in[n]$, $\bar{y}_i = u_i[\bar{x}_{-i}]^\circ(\bar{y}_i)$.\\ Let $d_i$ be a metric for $X_i$. From 
$y_{\boldsymbol{m}, i} = u[x_{\boldsymbol{m}, -i}]^\circ(y_{\boldsymbol{m}, i})$ we can write

\medskip\noindent
$\left\{
\begin{array}{lll}
 d_i\big(\bar{y}_i,  u_i[\bar{x}_{-i}]^\circ(\bar{y}_i)\big) \leqslant  
 d_i\big(\bar{y}_i, y_{\boldsymbol{m}, i}\big)\quad + \quad d_i\big(u[x_{\boldsymbol{m}, -i}]^\circ(y_{\boldsymbol{m}, i}), u[x_{\boldsymbol{m}, -i}]^\circ(\bar{y}_ {i}) \big)\quad + &   &   \\
 & &\\
 \hspace{6.1cm}+\quad d_i\big(u[x_{\boldsymbol{m}, -i}]^\circ(\bar{y}_ {i}), u[\bar{x}_{-i}]^\circ(\bar{y}_i)\big).& &
\end{array}\right.
$

\medskip\noindent
Since $X_i$ is compact and $u[x_{\boldsymbol{m}, -i}]^\circ$ is continuous  we can, for $\varepsilon > 0$ choose $\delta_{\boldsymbol{m}}(\varepsilon) > 0$ such that, if 
$d_i(z_i, w_i) < \delta_{\boldsymbol{m}}(\varepsilon)$ then  
$d_i\big(u[x_{\boldsymbol{m}, -i}]^\circ(z_{ i}), u[x_{\boldsymbol{m}, -i}]^\circ(w_ {i}) \big) < \varepsilon/3$. 

\medskip\noindent We can choose $\boldsymbol{m} = \boldsymbol{m}(\varepsilon)$ such that 
$d_i\big(\bar{y}_i, y_{\boldsymbol{m(\varepsilon)}, i}\big)< \varepsilon/3$ and\\ 
$d_i\big(u[x_{\boldsymbol{m(\varepsilon)}, -i}]^\circ(\bar{y}_ {i}), u[\bar{x}_{-i}]^\circ(\bar{y}_i)\big)< \varepsilon/3$.

\medskip\noindent
We have shown that $ d_i\big(\bar{y}_i,  u_i[\bar{x}_{-i}]^\circ(\bar{y}_i)\big) \leqslant  \varepsilon$. \hfill$\Box$

\begin{lem}\label{homotoptrivE} Assume that the strategy spaces are all compact  inf-semilattices with path-connected intervals. If, for all $x\in\prod_{i\in[n]}X_i$ and all $i\in [n]$,\\ $u_i[x_{-i}]^\circ : X_i\to X_i$ is continuous then, for all $x\in\prod_{i\in[n]}X_i$, the set $\mathbb{E}(x)$ is homotopically trivial.
\end{lem}
{\it Proof:} For all $x\in\prod_{i\in[n]}X_i$, $\mathbb{E}(x)$ is a compact inf-semilattice; we have seen that it has a smallest element and that it is path-connected. By a theorem of D.R. Brown, Theorem B in \cite{brown}, 
$\mathbb{E}(x)$ is homotopically trivial. \hfill$\Box$

\bigskip\noindent
A topological inf-semilattice is a {\bf Lawson semilattice} if each point has a neighbourhood base consisting of inf-semilattices. For example, a sub-inf-semilattice of $\mathbb{R}^n$ is a Lawson semilattice. A product of Lawson semilattices equipped with the product topology is a Lawson semilattice.

\begin{thm}\label{existnasheff} Let $\mathcal{G} = \big(X_i, u_i\big)_{i\in [n]}$ be a quasi-Leontief game such that:

\medskip\noindent $(1)$ all the strategy spaces are metrizable compact and locally connected path connected Lawson semilattices; 

\medskip\noindent $(2)$ all the payoff functions $u_i: \prod_{i\in [n]}X_i\to\mathbb{R}$ are continuous; 

\medskip\noindent $(3)$  for all $x\in\prod_{i\in[n]}X_i$ and all $i\in [n]$,  $u_i[x_{-i}]^\circ : X_i\to X_i$ is continuous. 

\medskip\noindent Then, $\mathcal{G}$ has an efficient Nash point.

\end{thm}

{\it Proof:} From the preceeding lemmas, $\mathbb{E}$ is an upper semicontinuous map with non empty homotopically trivial values from 
$\prod_{i\in [n]}X_i$ to itself. \\
Each $X_i$ has path connected intervals therefore, given two arbitrary points $x_1$ and $x_2$ of $X_i$ there is a continuous path from $x_1$ to $x_1\wedge x_2$ and a continuous path from $x_1\wedge x_2$ to $x_2$; $X_i$ is path connected, and {\it a fortiori} connected. By a theorem of M. McWaters, \cite{macwa}, condition $(1)$ implies that each $X_i$ is an absolute retract; $\prod_{i\in [n]}X_i$ is therefore an absolute retract. \\
Homotopically trivial sets being acyclic,  $\mathbb{E}$ is an upper semicontinous map with non empty acyclic values from the compact absolute retract $\prod_{i\in [n]}X_i$ to itself. By the Eilenberg-Montgomery 
Theorem, \cite{dugra} Page 543 Corollary (7.5), there exists $x\in\prod_{i\in [n]}X_i$ such that $x\in \mathbb{E}(x)$.\hfill$\Box$ 

\begin{cor}\label{nasheffn}
Let $\mathcal{G} = \big(X_i, u_i\big)_{i\in [n]}$ be a quasi-Leontief game such that:

\medskip\noindent $(1)$ $X_i$ is a  compact  inf-convex subset of $\mathbb{R}^{n_i}$; 

\medskip\noindent $(2)$ all the payoff functions $u_i: \prod_{i\in [n]}X_i\to\mathbb{R}$ are continuous; 

\medskip\noindent $(3)$  for all $x\in\prod_{i\in[n]}X_i$ and all $i\in [n]$,  $u_i[x_{-i}]^\circ : X_i\to X_i$ is continuous. 

\medskip\noindent Then, $\mathcal{G}$ has an efficient Nash point.
\end{cor}

{\it Proof:} In $\mathbb{R}^n$ intervals are path connected; an inf-convex subset of $\mathbb{R}^n$ is therefore a subsemilattice with path connected intervals. Let $x$ be an arbitrary point of $X_i$ and let $U$ be a  neighbourhood of $x$ in $X_i$. Choose a neighbourhood $V$ of $x$ in  
$\mathbb{R}^{n_i}$ such that $U = X_i\cap V$ and a neighbourhood $W$ 
of $x$ in  
$\mathbb{R}^{n_i}$ such that $W\subset V$ and $W$ is inf-convex; $W$ could for example be a box around $x$. Since the intersection of two inf-convex sets is inf-convex and an inf-convex set is path connected we have that $W\cap X_i$ is a neighbourhood of $x$ in $X_i$ that is a connected  inf-semilattice contained in $U$; this shows that $X_i$ is a locally connected Lawson semilattice. \hfill$\Box$ 

\bigskip\noindent {\bf A final remark:} The algebraic and the topological assumptions used throughout this paper are not as different as one could believe. The compactness assumption is natural and at the same time seems to be somewhat indeterminate but, as a matter of fact, there is at most one topology on a given semillatice for which it is a compact topological semilattice and continuity is defined entirely in terms of order convergence. For a lattice, that topology is explicitely determined by the algebraic structure of the lattice.  For this, and more and topological lattices and  semilattices see Theorem 15 and Corollary 16 in \cite{law}.

\end{document}